\documentclass[10pt]{amsart}
\usepackage{latexsym}
\usepackage{mathrsfs,amsmath,amssymb,amscd,amsthm,amsfonts,verbatim,eufrak}

\usepackage[dvips]{epsfig}
\usepackage{color}

\def\1{{\hat{1}}}
\def\0{{\hat{0}}}

\def\vees{{\;\vee^{\scriptscriptstyle 2 }\;}}
\def\veek{{\;\vee^{\scriptscriptstyle k }\;}}

\def\Tk{{\mathcal{T}^{k}_n}}
\def\P{{\Pi_m}}
\def\Pk{{\Pi^{(k)}_m}}

\def\FF{{\mathscr F}}
\def\GG{{\mathcal G}}
\def\HH{{\mathcal H}}
\def\II{{\mathcal I}}
\def\LL{{\mathcal L}}

\def\NN{{\mathcal N}}

\def\TT{{\mathcal T}}

\def\0{{\hat{0}}}

\newtheorem{thm}{Theorem}[section]
\newtheorem{df} [thm]{Definition}
\newtheorem{lm}  [thm]{Lemma}
\newtheorem{crl} [thm]{Corollary}
\newtheorem{prop}[thm]{Proposition}

\newtheorem{rem}[thm]{Remark}

\newtheorem{expl}[thm]{Example}
\newtheorem{cstr}[thm]{Construction}
\numberwithin{equation}{section}

\newenvironment{pf}{\noindent {\bf Proof. }}{\hfill $\Box$\vspace{0.3cm}}

\begin{document}

\title[Subdivision of complexes of $k$-Trees]{Subdivision of complexes of $k$-Trees}

\author{Emanuele Delucchi}

\address{Department of Mathematics, ETH Zurich, 8092 Zurich, Switzerland}
\email{delucchi@math.ethz.ch}

\thanks{{MSC 2000 Classification:} primary 05E25, secondary 57Q05.}

\begin{abstract}
Let $\Pi^{(k)}_{(n-1)k+1}$ be the poset of partitions of $\{1,2,\dots,(n-1)k+1\}$ with block sizes congruent to $1$ modulo $k$. We prove that the order complex $\Delta(\Pi^{(k)}_{(n-1)k+1})$ is a subdivision of the complex of $k$-trees $\Tk$, thereby answering a question posed by Feichtner \cite[5.2]{F}.

The result is obtained by an {\em ad-hoc} generalization of concepts from the theory of nested set complexes to non-lattices.
\end{abstract}

\maketitle


\section*{Introduction}
\label{sect_intr}


The complex of $k$-trees $\Tk$ is the abstract simplicial complex whose faces correspond to combinatorial types of rooted trees with $m=(n-1)k+1$ labelled leaves and all outdegrees congruent to $1$ modulo $k$ (and larger than $k$). Face relations among these trees are such that $T_2$ is a face of $T_1$ if and only if $T_2$ is obtained from $T_1$ by contraction of internal edges. Complexes of $k$-trees are a generalization of complexes of trees \cite{Bo} defined by Hanlon in \cite{H}. The complexes $\Tk$ are pure simplicial complexes of dimension $n-3$. They were shown to be Cohen-Macaulay by Hanlon in \cite{H}, and shellability was proven independently by Trappmann and Ziegler \cite{TZ} and Wachs \cite{W}.

The symmetric group $S_m$ acts naturally on $\Tk$ by permutation of the leaves. In fact, the complex $\Tk$ was introduced by Hanlon because of the relation of its top homology group 
to the $1^m$ homogeneous component of the free Lie-$k$-algebra with $m$ odd  generators, both as $S_m$-modules. Lie-$k$-algebras are a generalization of Lie algebras introduced by Hanlon and Wachs in \cite{HW}. Here, the poset $\Pk$ of partitions of $\{1,\dots , m\}$ with block sizes congruent to $1$ modulo $k$ enters our picture. 
This poset was shown to be Cohen-Macaulay by Bj\"orner \cite{Bj}, and was studied by Calderbank, Hanlon and Robinson in \cite{CHR} from the point of view of representation theory.
Combining the results of \cite{H} and \cite{HW} one finds that there is a $S_m$-module isomorphism between $\tilde{H}_{n-3}(\Tk)$ and $\tilde{H}_{n-3}(\Pk)$, since both are isomorphic to the $1^m$ homogeneous piece of the free Lie-$k$-algebra with $m$ odd generators.

This situation leads to the very natural question of the geometric relationship between the two complexes.

The case $k=1$ is interesting on its own right. The complex $\TT_n=\TT_n^1$ (called {\em complex of trees})  first appeared in \cite{B}, in the context of homotopy theory. 
Its $S_m$-representation theory was studied by Robinson and Whitehouse \cite{RW}, it figured in geometric group theory in work of Vogtmann \cite{V}, and recently attracted interest again in the study of phylogenetic trees \cite{BHV}. In fact,  the order complex of the partition lattice $\Pi_m$ is PL-homeomorphic to $\Tk$ \cite{AK}. More specifically, $\Delta(\Pi_m)$ is obtained from $\TT_n$ by a sequence of stellar subdivisions \cite{F}. 
This is shown by applying the theory of {\em combinatorial nested set complexes}, introduced by Feichtner and Kozlov in \cite{FK} as the combinatorial framework and counterpart of De Concini - Procesi models of subspace arrangements \cite{DP}. In \cite[Question 5.2]{F}, Feichtner asks whether for general $k$ the complexes $\Pk$ and $\Tk$ are related by a subdivision. 

In this note we answer this question positively. The direct approach via nested set complexes does not work: other than $\P$, $\Pk$ is not a lattice, which is an essential condition for the applicability of the theory of nested sets. However, in our context we can introduce analogous concepts that allow us to reduce the problem to a situation that can be handled by the existing theory.

In Section \ref{nsc} we review the definition and some of the basic properties of nested set complexes.
In Section \ref{combsub} we set up the topological framework for subdivisions of simplicial complexes and, with Lemma \ref{lemma}, state the key tool that will enable us to prove our result.   
In Section \ref{terza} the analogies with the classical framework of nested sets will be put to work, giving explicite subdivisions of each simplex of $\Tk$. By Lemma \ref{lemma} we will then be able to construct a subdivision of the whole complex. Although our method is very much tailored to the specific situation of the partition posets, we have kept the notations and the arguments as close as possible to the original theory of nested sets, hoping that this may be a first step towards a generalization of the theory to non-semilattices.

\section{Nested set complexes}
\label{nsc}
We recall the basic definitions and some useful facts about nested set complexes. For a more extensive account see \cite{FK}.

{\df \label{building} Let $\LL$ be a meet-semilattice. A {\bf building set} of $\LL$ is a subset $\GG\subseteq\LL\setminus\{\hat{0}\}$ such that for any $x\in\LL\setminus\hat{0}$ there is an isomorphism
\begin{displaymath}
\varphi_x: {\Large\prod_{\scriptstyle{{z\in\textrm{\textup{max}}}\,\GG_{\leq x}}}} [\hat{0},z] \rightarrow [\hat{0},x]
\end{displaymath}
with $\varphi_x(\0,\dots\0,z,\0,\dots,\0)=z$ for $z\in\textrm{\textup{max}}\,\GG_{\leq x}$.}

\vspace{4mm}

Given a building set $\GG$, the set $F_\GG(x):=\textrm{\textup{max}}\,\GG_{\leq x}$ is called the set of {\em factors of $x$ in $\GG$}. We may sometimes drop the suffix and write only $F(x)$ if no confusion can arise.\\

\noindent{\bf Remarks:}
\begin{itemize}
\item[(1)]Building sets can be defined as in \cite[2.3, (2)]{FK}, independently of the existence of a join operation, for any poset with a unique minimal element. 
\item[(2)]If a join operation exists in $\LL$, the isomorphism $\varphi_x$ can always be chosen such that $\varphi_x(z_1,\dots,z_n)=z_1\vee\dots\vee z_n$ \cite[Proof of 2.3]{FK}. 
\end{itemize}

{\df\label{nested} Let $\LL$ be a semilattice and $\GG$ a building set in $\LL$. We call a set $N\subseteq\GG$ {\bf nested} ($\GG$-nested, if specification is needed) if, for any set $\{x_1,\dots,x_\ell\}\subseteq N$ ($\ell\geq 2$) of incomparable elements, the join $x_1\vee \dots\vee x_\ell$ exists and is not an element of $\GG$.
The {\bf nested set complex} of $\LL$ with respect to $\GG$, denoted $\NN(\LL,\GG)$, is the abstract simplicial complex of all nonempty $\GG$-nested sets (with the apex removed, if $\LL$ has a maximal element).
}

\vspace{4mm}

\noindent{\bf Remarks:}
\begin{itemize}
\item[(1)]This definition implies that $F(x_1\vee\dots\vee x_\ell)=\{x_1,\dots,x_\ell\}$ for any set of incomparable elements $\{x_1,\dots,x_\ell\}\subseteq N$, see \cite[2.8, (2)]{FK}.
\item[(2)]Consider the maximal building set $\GG=\LL$. In this case, a set is $\GG$-nested if and only if it is linearly ordered. Therefore we have $\NN(\LL,\LL)=\Delta(\LL)$, the order complex of (the proper part of) $\LL$. 
\end{itemize}

If $\LL$ is atomic, then Feichtner and Yuzvinsky \cite{FY} give an explicit geometric realization of nested set complexes of $\LL$. If $\HH\subseteq\GG$ are building sets of $\LL$, Feichtner and M\"uller \cite{FM} prove that $\NN(\LL,\GG)$ arises from $\NN(\LL,\HH)$ by a sequence of stellar subdivisions.

\section{Subdivisions of abstract simplicial complexes}
\label{combsub}


Let $K$ be an abstract simplicial complex on the vertex set $V(K)$. We denote its face poset by $\FF(K)$.
Given $s\in\FF(K)$, we will write $K_s$ for the simplicial complex of all nonempty subsets of the vertex set $V(K_s)=s$.

Recall that the \emph{realization} of a simplicial complex $K$ is the topological space $|K|$ whose points are the formal convex combinations $\Sigma t_v v$ of vertices (i.e., the coefficients $t_v$ are all nonnegative and $\sum_{v\in V(K)} t_v\! =\! 1$) such that the set $\{v \in V(K) |\, t_v \neq 0\}$ belongs to the family $K$. For any simplex $s\in K$ define $|s|:=|K_s|=\{\Sigma t_v v \in |K| |\, t_v \neq 0 \Rightarrow v\in s \}$. If $s\in K$ contains $n$ vertices there is a natural correspondence between $|s|$ and the standard Euclidean $n$-simplex $\Delta^n$. The topology on $|K|$ is such that $A\subseteq |K|$ is closed if and only if $A\cap |s|$ is closed (in the Euclidean topology induced by $\Delta^n$) for each simplex $s$ of $K$.

The interior of a simplex $s\in K$ is then $\langle s \rangle := |s|\setminus\partial |s|$, and for each point $x\in|K|$ we can now define the {\em carrier} of $x$ as $\tau(x):=\textrm{min}\{s\in K | x\in \langle s \rangle\}$.\\




{\df{\em (see \cite[p.121]{Sp}).} Let $P,Q$ be simplicial complexes. $P$ is a {\bf subdivision} of $Q$ if there exists a poset map $\phi : \FF(P) \rightarrow \FF(Q)$ and an (injective) mapping $f^0: V(P)\rightarrow |Q|$ such that $f^0(V(P_p))\subseteq |Q_{\phi(p)}|$ for all simplices $p\in P$, and the linear extension $f:|P|\rightarrow |Q|$ of $f^0$ to $|P|$ is a homeomorphism.
}\\

The map $\phi$ will be called the \emph{carrier map}, and $f^0$ the \emph{vertex map} of the subdivision.

%

Our point of view will be to take a given map between the face posets, and look for conditions under which this map is the carrier map of a subdivision (or, as we will say in the following: \emph{is} a carrier map).

Suppose that we have a good candidate for $f^0$. The linear extension $f$ is continuous because it is linear -hence continuous- on every closed simplex. Realizations of simplicial complexes are compact, Hausdorff spaces. Therefore it suffices to check bijectivity of $f$ to conclude that $f$ is a homeomorphism. We summarize this in the following lemma for later reference. 

{\lm{\em (see \cite[p.122]{Sp})} Let $P,Q$ be simplicial complexes, and let a vertex map $f^0:V(P)\rightarrow |Q|$ be given. The linear extension 
\begin{equation*}\begin{array}{rccc}
f:& |P| &\rightarrow & |Q| \\
  & \Sigma t_vv & \mapsto & \Sigma t_vf^0(v) 
\end{array}\end{equation*} 
of $f^0$ to $|P|$ is a homeomorphism if and only if for every simplex $q$ of $Q$ the set $\{f(\langle p \rangle) | f^0(p)\subseteq |q|\}$ is a finite partition of $\langle q \rangle$.}\\


We now define the concept of {\em compatibility} of carrier maps. This will lead to a lemma that will be one of the key tools in our argument.

{\df \label{compatible carrier map} Let $P,Q$ be simplicial complexes, $P_1,P_2$ and $Q_1,Q_2$ subcomplexes of $P$, resp. of $Q$. Two poset maps $\phi_1:\FF(P_1)\rightarrow \FF(Q_1)$, $\phi_2:\FF(P_2)\rightarrow \FF(Q_2)$ are called {\bf compatible carrier maps} if for $i=1,2$ there are vertex maps $f^0_i:V(P_i)\rightarrow V(Q_i)$ so that $P_i$ is a subdivision of $Q_i$ with carrier map $\phi_i$, and moreover $f^0_1 = f^0_2$ on $V(P_1)\cap V(P_2)$.}

{\lm \label{lemma} A poset map $\phi:\FF(P)\rightarrow\FF(Q)$ is the carrier map of a subdivision if the maps $\phi_q : \phi^{-1}(\FF(Q)_{\leq q})\rightarrow \FF(Q)_{\leq q}$ are compatible carrier maps for all $q \in \FF(Q)$.
}

\vspace{2mm}
{\pf Compatibility of the $\phi_q$ ensures that we can choose the $f^0_q$ such that if $q_1\leq q_2$, $f^0_{q_1}(v)=f^0_{q_2}(v)$ for all $v$ where it is defined. Then, for all $v\in V(P)$, we define
\begin{center}
$f^0(v):=f^0_{\phi(v)}(v).$
\end{center}
Clearly $f^0(p)=f^0_{\phi(p)}(p)\subseteq |\phi(p)|$ for all simplices $p\in P$. Now take any $q\in Q$ and note that in particular $f=f_q$ on $\{v\in V(P)|\, \{v\}\in \phi^{-1}(q)\}$. Hence $\{f(\langle p \rangle) |\, f(p)\subseteq |q|\}=\{f_q(\langle p\rangle)| \, p\leq\phi_q(q)\}$, and this is by assumption a finite partition of $|q|$. \hfill$\square$
}\\

The condition of compatibility may be hard to check, and therefore Lemma \ref{lemma} does not give a handy tool in general. However, it becomes useful in our situation, where compatibility will emerge from an explicit description of 'partial' subdivision processes.

\section{$\Delta(\Pk)$ is a subdivision of $\Tk$}
\label{terza}

\subsection{Complexes of $k$-trees.} We recall the definition of complexes of $k$-trees from \cite{H} and view them in the context of a nested-set-like construction.

{\df \hspace{-4pt}{\bf([H])} For positive integers $k,n$ let $\Tk$ denote the collection of all trees having $(n-1)k+2$ labelled leaves and having the property that every internal vertex has degree $rk+2$ for some $r>0$. We can partially order this set by saying that $T_1 > T_2$ if $T_2$ can be obtained from $T_1$ by contraction of some set of internal edges. }

\begin{figure}[htbp]
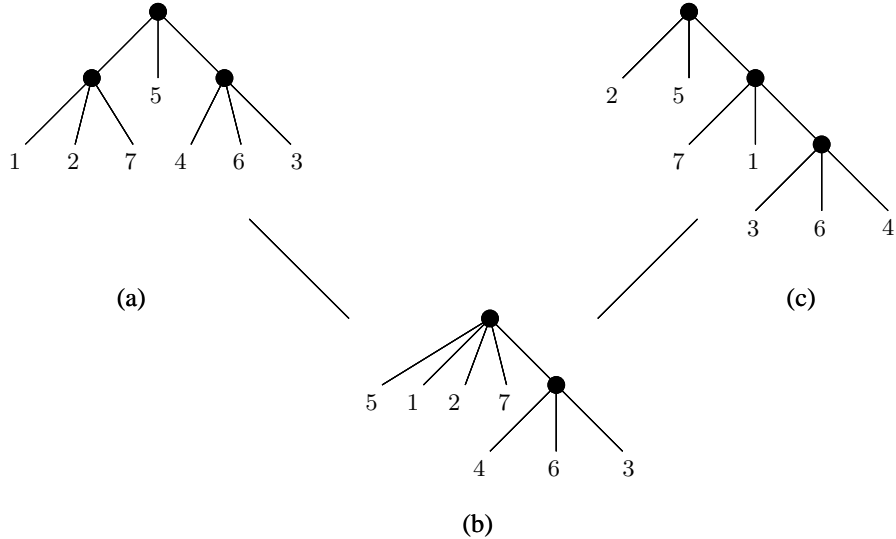

\begin{center}
  \begin{picture}(0,0)%
    \includegraphics{ktree.pstex}%
  \end{picture}%
  \input{ktree.pstex_t}%
   
\caption{Two maximal cells of $\TT_4^2$ ((a),(c)) and their common point ((b)).}
\end{center}
\end{figure}

\vspace{2mm}

For us it will be convenient to work with simplicial complexes and rooted trees. We rewrite the previous definition, and let $\Tk$ denote the abstract simplicial complex of rooted trees on $m=(n-1)k+1$ leaves with all outdegrees larger than $1$ and congruent to $1$ modulo $k$. The face relations are such that $T_2$ is a face of $T_1$ if $T_2$ is obtained from $T_1$ by a contraction of a set of internal edges.

\vspace{4mm}

Let $\II$ denote the minimal building set of $\P$. It consists of all partitions with only one block of size bigger than $1$.

At this point it is worth to remind that $\Pk$ is not a lattice. For example, the join of the partitions $(123)4567$ and $1(234)567$, considered as elements of $\Pi_7$, is clearly $(1234)567$. But as elements of $\Pi_7^{(2)}$ they have three different minimal upper bounds: $(12345)67$, $(12346)57$ and $(12347)56$.

This is the reason why the theory of combinatorial nested sets cannot be directly applied here. Nevertheless, in the following we will reduce our problem, at least 'locally', to a situation where the existing theory can be applied.

We define a subset of $\Pk$,
\begin{center}$\GG:=\{x\in\II \vert\; rk(x)\equiv 0 \pmod{k} \} \subseteq\Pk$.\end{center} 

It will serve as an analogue of the building set $\II$ in our non-lattice context.

{\prop There is a bijective correspondence between the simplicial complex $\Tk$ and  the family of sets $\{N\in \NN(\P,\II)\vert N\subseteq\GG\}$. In particular, the latter is a simplicial complex.}

\vspace{2mm}

{\pf  We follow \cite[proof of 3.1]{F}, where an explicit correspondence is given in the case of $k=1$. \\
\indent ``$\subseteq$'':  Let $T\in\Tk$ be a (rooted) tree with inner vertices $t_1,\dots,t_k$. The set of leaves below the vertex $t$ will be denoted by $\ell(t)$. For $N(T):=\{\ell(t_i) \vert i=1,\dots,k\}$ we have $N(T)\in \NN(\P,\II)$ \cite{F}, and, by the condition on outdegrees, clearly $N(T)\in\GG$.\\
\indent ``$\supseteq$'': To every $N=\{S_1,\dots,S_k\}\in\NN(\P,\II)$ we can construct a tree $T(N)$ as in \cite{F}. The set of internal vertices is $\tilde{N}=N\cup\{[m]\}$ (where $[m]$ corresponds to the root), and the covering relation is defined by $U_1<U_2$ if and only if $U_2\in\textrm{max} N_{<U_1}$ for all $U_1,U_2\in\tilde{N}$. The outdegree $out(U)$ of an internal vertex $U\in \tilde{N}$ is:
\begin{center}$out(U)=|U\setminus\bigcup_{S\in N, S<U} S|+|\{S\in N \vert S<U\}|=|U|-\sum_{S\in N, S<U} |S| +|\{S\in N \vert S<U\}|$. \end{center} 
Moreover, if $N\in\GG$, then $|S|\equiv 1 \pmod{k}$, and with $s:=|\{S\in N \vert S<U\}|$ we have $out(U)\equiv 1 -s +s \equiv 1 \pmod{k}$. \hfill $\square$ }

\vspace{4mm}

In the following we will consider $\Pk$ as a subset of $\P$. Given two elements $a,b \in \Pk$ we shall denote by $a\vee b$ the join of them as elements of $\P$, and by $a \veek  b$ we will mean the set of all minimal upper bounds of $a$ and $b$ in $\Pk$. In the example above we would then have $(123)4567 \vee 1(234)567 = (1234)567$, $(123)4567\vees 1(234)567 = \{(12345)67, (12346)57, (12347)56\}$. If it happens that $a$ and $b$ have a unique least upper bound we will write $\exists ! \; a\veek b$, and then denote by $a \veek b$ the (unique) element $x\in a\veek b$.

We write $F(x)$ for $F_{\II}(x)$, the set of factors of $x$ with respect to the building set $\II$ of $\P$ (see Definition \ref{building}).

{\lm {\bf\emph{(Properties of $\GG \subseteq \Pk$)}}\label{G}
\begin{itemize}
\item[a)] For all $x\in\Pk$, we have $F(x)=\textrm{max}\,\GG_{\leq x}$.
\item[b)] Let $A,B$ denote the non-singleton blocks of two partitions $a,b\in\GG$. Then we have:
\begin{center}$A\cap B = \emptyset \;\; \Leftrightarrow \;\;\exists ! \; a\veek b \;\textrm{and}\; a\veek b \notin \GG$. 
\end{center}
\end{itemize}
}

\vspace{2mm}\pagebreak
{\pf 
\begin{itemize}
\item[a)]The only non-singleton block of any partition in $F(x)$ is always a block of $x$, therefore has size equal to $1$ modulo $k$.
\item[b)]If $A\cap B = \emptyset$, the only non-singleton blocks of $a\vee b$ are $A$ and $B$, and therefore $a \vee b \in \Pk$. So $\exists ! \;a\veek b$, and since $a\vee b \notin \II$ we have $a \veek b \notin \GG$.

\noindent Suppose $A\cap B \neq \emptyset$. Then $a\vee b$ has only one non-singleton block, namely $A\cup B$. If $\vert A \cup B\vert \equiv 1 \pmod{k}$ then $\exists ! \; a \veek b$, but $a\veek b \in \GG$. If $\vert A \cup B \vert \not\equiv 1  \pmod{k}$, then $\vert a \veek b \vert >  1$. In both cases the claim follows.\hfill $\square$
\end{itemize}
}
\vspace{4mm}

This Lemma suggests the following definition, which mimics Definition \ref{nested}. 

{\df Let $\NN$ denote the family of all subsets $N\subseteq \GG$ such that for any set of more than two incomparable elements $ \{x_1,\dots , x_\ell\} \subseteq N$ we have that $\exists ! \; x_1\veek\dots\veek x_\ell$ and $ x_1\veek\dots\veek x_\ell \notin \GG $.}

\vspace{4mm}

The following statement can be viewed as the analogue of \cite[Theorem 3.1]{F} in our context. After the above preparations, the proof is an easy check of the definitions.

{\lm With the identification of the elements of $\GG$ with the vertices of $\Tk$ given above, $\NN = \Tk$ as simplicial complexes. \hfill $\square$\label{TN} }

\vspace{4mm}

We carry on the analogy. For $x\in\Pk$ let us now write $F^k(x):=\textrm{max}\;\GG_{\leq x}$, and given a chain $\omega\in\Delta(\Pk)$, let $F^k(\omega):=\bigcup_{x\in\omega} F^k(x)$. Then we have following suggestive corollary:

{\crl For any $\omega \in \Delta(\Pk)$, $F^k(\omega)\in \NN$.}

\vspace{2mm}
{\pf First note that $F^k(\omega)\in\Tk$, because $F^k(\omega)\subseteq \GG $ by Lemma \ref{G}, and $F^k(\omega)\in \NN(\P , \II)$ by \cite{FK}[2.8 (3)]. Then apply Lemma \ref{TN}.\hfill $\square$}

\subsection{The poset $\Sigma(N)$} For $N \in \NN$ consider 
\begin{equation*}
\Sigma(N):=\Big\{\bigvee_{x\in  X}\!\!{}^{\stackrel{\textstyle k}{}}\, x \vert \; X\subseteq N \Big\} \subseteq \Pk,
\end{equation*}
the set of all joins of subsets of $N$, with the induced order (the join over the empty set gives the minimal element).
Note that $N\in\NN$ implies in particular that 
for every $A \subseteq N$ the join $\bigvee^k_{x\in A}x$ exists and is unique. The join in $\Sigma(N)$ will be denoted by $\vee$.

{\rem \label{rmk} For all $a\in\Sigma(N)$ we have $a=\bigvee\!{}^{\stackrel{\scriptstyle k}{}} N_{\leq a}=\bigvee\textrm{max}\,N_{\leq a}$.}


{\prop 
The poset $\Sigma(N)$ is a lattice, and $N$ is a building set for $\Sigma(N)$.
}

\vspace{2mm}
{\pf 
We divide the proof in two parts.
\begin{description}
\item[Lattice] $\Sigma(N)$ has a maximal and a minimal element (the join over all $N$ and the join over the empty set). For the meet operation take $a,b\in \Sigma(N)$ and define 
\begin{center}$A:= \textrm{max}\, (N_{\leq a}\cap N_{\leq b})$.\end{center}
This is a set of incomparable elements in $N$, so $\exists ! \;\bigvee^k A$. We claim that this is the minimal lower bound for $a$ and $b$.\\
Indeed, for every $c\in\Sigma(N)$ with $c\leq a$ and $c\leq b$, clearly $N_{\leq c }\subseteq N_{\leq a}$ and  $N_{\leq c}\subseteq N_{\leq b}\,$. Therefore $c=\bigvee^k\textrm{max}\,  N_{\leq c} \leq \bigvee^k A$. 
\vspace{2mm}
\item[N is building set]Take any $x\in\Sigma(N)$ and let $\{x_1, \dots
  , x_\ell\}=\textrm{max}\, N_{\leq x}$. The map
\begin{equation*}\begin{array}{rccc}
\phi :&\prod_{j=1}^{\ell} [\hat{0}, x_j]& \rightarrow &[\hat{0}, x]\\
&(a_1,\dots , a_\ell)&\mapsto & a_1 \vee \dots \vee a_\ell
\end{array}\end{equation*}
is surjective by definition of $\Sigma(N)$.

For injectivity, note that $\{x_1, \dots, x_l\}\in \NN$ as a set of incomparable elements, and by \ref{G} we then know the
block decomposition of the partition $x=\bigvee^k x_j=\bigvee
x_j$. From this we conclude $F(x)=\{x_1,\dots , x_l\}$, and the
defining property of building sets in lattices gives an isomorphism $\psi
: \prod [\hat{0}, x_j] \rightarrow [\hat{0}, x]$ \emph{in the full
  partition lattice}, defined
by $\psi(x_1,\dots , x_\ell)=x_1 \vee  \dots \vee x_\ell$ (see \cite{FK}).

Now suppose $\phi(a)=\phi(b)$ for $a\neq b$. By construction then $\psi(a)=\phi(a)=\phi(b)=\psi(b)$ - which contradicts bijectivity of $\psi$. \hfill $\square$


\end{description}
}
\vspace{2mm}
Recall that, given a face $s$ of a simplicial complex $K$, we denote by $K_s$ the subcomplex of all faces of $s$.

{\crl\label{stellar} $\Delta (\Sigma (N))$ is obtained from $\NN_{N}$ by a
  sequence of stellar subdivisions.}
\vspace{2mm}

{\pf This is obtained as the main result of \cite{FM}, because $\NN_{ N} =
  \NN(\Sigma(N), N)$. The atomicity condition - needed in \cite{FM} to construct an explicit realization of the complex - can here be dropped: $\NN_{N}$ is a simplex and we have a standard geometrical realization for it. The subdivisions can be still described as combinatorial blowups in the face poset. \hfill $\square$}

\vspace{4mm}

We have to take a closer look at the carrier maps of that subdivision.

{\prop\label{blowup} Let $\LL$ be an atomic semi-lattice, $\HH\subseteq\GG$ two building sets of $\LL$. The carrier $\tau(N)$ of a simplex $N$ of the nested set complex $\NN(\LL,\GG)$ in the subdivision of $\NN(\LL,\HH)$ by $\NN(\LL,\GG)$ is given by
\begin{displaymath}
\tau(N)=\bigvee_{x\in N} F_{\HH}(x) \; \in \FF(\NN(\LL,\HH)).
\end{displaymath}}

\vspace{2mm}

{\pf Every linear extension $h_1,h_2,h_3,\dots, h_\ell$ of the order of $\LL$ on $\GG\setminus\HH$ gives a sequence of stellar subdivisions leading from $\NN(\LL,\HH)$ to $\NN(\LL,\GG)$ (see \cite{FM}). We proceed by induction on $k$.

If $\ell=1$, then $\NN(\LL,\GG)$ is obtained from $\NN(\LL,\HH)$ by stellar subdivision of the simplex $F_{\HH}(h_1)$, so that the vertex map $f^0: V(\NN(\LL,\GG))\rightarrow |\NN(\LL,\HH)|$ maps $h_1$ to the barycenter of $F_{\HH}(h_1)$, whereas $f^0(v)=|v|$ for all $v\neq h_1$. The carrier of $N\in\NN(\LL,\GG)$ is the smallest simplex $\tau(N)$ in $\NN(\LL,\HH)$ such that $f^0(N)\subseteq |\tau(N)|$. This is the same to say that $\tau(N)$ is the smallest simplex containing the carrier of every vertex of $N$, hence $\tau(N)=\bigvee_{x\in N} F_{\HH}(x)$.

Now let $\ell>1$, define $\HH_{\ell}:=\GG\setminus\{h_\ell\}$ and let $\tau_\ell(N)$ be the carrier of $N\in\NN(\LL,\GG)$ in $\NN(\LL,\HH_\ell)$. 
By the argument given in the case $\ell=1$, we have $\tau_\ell(N)=\bigvee_{x\in N} F_{\HH_\ell}(x)$. 

Note that, with respect to $\NN(\LL,\HH)$, the carrier of $N\in\NN(\LL,\GG)$ is the same as the carrier of $\tau_\ell(N)\in\NN(\LL,\HH_\ell)$. Using the induction hypothesis we can now write:

\begin{displaymath}
\tau(N)=\bigvee_{y\in \tau_\ell(N)} F_\HH (y) = 
\left\{\begin{array}{ll}
\bigvee_{x\in N} F_\HH (x) & \textrm{if } h_\ell \not\in N, \\
\bigvee_{x\in N\setminus\{h_\ell\}} F_\HH(x) \vee F_\HH(F_{\HH_\ell}(h_\ell)) & \textrm{if } h_\ell \in N.
\end{array}\right.
\end{displaymath} 

Because $\HH \subset \HH_\ell$, we have $ F_\HH(F_{\HH_\ell}(h_\ell))=F_{\HH}(h_\ell)$. Hence in any case the claim follows.$\square$
}
\vspace{4mm}

As in the proof of Corollary \ref{stellar}, the atomicity condition is needed in the general statement to construct an explicit geometric realization of the nested set complex. We again note that, in our case, the lattice $\Sigma(N)$ is not atomic but the nested set complexes to which we apply the result of \cite{FM} are always simplices. We have then following application of Proposition \ref{blowup}:

{\prop\label{carrmap} For $N\in\NN$, the map
\begin{equation*}\begin{array}{rccc}
\phi_N : & \FF(\Delta(\Sigma(N)))& \rightarrow& \FF(\NN_N) \\
         &    \omega     & \mapsto    &  F^{(k)}(\omega)
\end{array}\end{equation*}
is the carrier map of the subdivision of $\NN_N$ by $\Delta(\Sigma(N))$ described in Corollary \ref{stellar}.}

\vspace{2mm}

{\pf
Proposition \ref{blowup} gives $\phi_N(\omega)=\bigvee_{x\in\omega} F^{(k)}(x)$. So we only have to show that the right hand side equals $F^{(k)}(\omega)$. For this, it is enough to show

\begin{description}
\item[{\bf Claim}]{\em Let $\LL$ be a semi-lattice and $\GG$ a building set in $\LL$. Then $x>y\in\LL$ implies $F(x)\vee F(y)= F(x) \cup F(y)$ in $\FF(\NN(\LL,\GG))$. }
\item[Proof] We show that $F(x)\cup F(y)$ is nested. Since the set of factors of an element is nested, we only need to show that for two incomparable elements $z_1\in F(x)$ and $z_2\in F(y)$, the join $z_1\vee z_2$ is not in $\GG$. From $x>y$ we know that $z_1 \vee z_2 \leq x$ and there is exactly one $z\in F(x)$ such that $z_2\leq z$ (see \cite[2.5]{FK}). Recall the isomorphism $\varphi_x$ in definition \ref{building}: w.l.o.g. we can write $z_1\vee z_2 = \varphi_x(z_1,z_2,0,\dots,0)$, and this cannot be in $\GG$ because the preimage of any element of $\GG_{\leq x} = \bigcup_{z\in F(x)}[\hat{0},z]$ has only one nonzero entry.\hfill $\square$  
\end{description}
}

\subsection{The subdivision}After explicit study of the geometric structure of the subdivisions, we are now able to check the conditions of Lemma \ref{lemma}, and use it to prove our theorem.

 {\thm \label{teorema} $\Delta(\Pk)$ is a subdivision of  $\Tk$}
\vspace{2mm}

{\pf Consider the map 
\begin{equation*}\begin{array}{rccc}
\phi : & \FF(\Delta(\Pk))& \rightarrow& \FF(\NN) \\
         &    \omega     & \mapsto    &  F^{(k)}(\omega).
\end{array}\end{equation*}
It is easy to see that $\phi^{-1}(\FF(\NN)_{\leq N})=\FF(\Delta(\Sigma(N)))$, hence  by Proposition \ref{carrmap} the maps $\phi_N: \phi^{-1}(\FF(\NN)_{\leq N})\rightarrow \FF(\NN)_{\leq N}$ are carrier maps.

For compatibility note that $\Delta(\Sigma(N_1)) \cap \Delta(\Sigma(N_2)) = \Delta (\Sigma(N_1 \cap N_2))$ and that, for any $N$, $f^0_N (v)$ is the barycenter of $|F^{(k)}(v)|$. It is then clear that $f^0_{N_1}\!\!\upharpoonright_{V(\Delta(\Sigma(N_1))) \cap V(\Delta(\Sigma(N_2)))} = f^0_{N_1\cap N_2} = f^0_{N_2}\!\!\upharpoonright_{V(\Delta(\Sigma(N_1))) \cap V(\Delta(\Sigma(N_2)))}$.

With Lemma \ref{lemma} we conclude that $\phi$ is a carrier map.\hfill $\square$
}
\vspace{2mm}

{\rem $\;\!$} As in \cite{FM}, we see that every linear extension of $\Sigma(N)\setminus N$ gives a valid sequence of stellar subdivisions leading from $\NN_{N}$ to $\Delta(\Sigma(N))$. By the compatibility checked in \ref{carrmap} and \ref{teorema} it is easy to see that a valid sequence of stellar subdivisions leading from $\Tk$ to $\Delta(\Pk)$ can be obtained from any linear extension of $\Pk\setminus \GG$.   

\vspace{4mm}

To conclude this note let us consider the natural action of the symmetric group $S_m$ on the complex $\Tk$ (by permutation of leaves) and on the poset $\Pk$ (by permutation of elements). The correspondence of Lemma \ref{TN} is obviously $S_m$-equivariant. One can check that the carrier maps of the (whole) subdivision are also $S_m$-equivariant, even if at some single subdivision step this may not be true. We then recover the result of \cite{H}, \cite{HW}:

{\crl The top degree homology modules of $\Tk$ and $\Pk$ are isomorphic as $S_m$-modules:
\begin{displaymath}
\widetilde{H}_{n-3}(\Tk) \cong_{S_m} \widetilde{H}_{n-3}(\Pk).
\end{displaymath}}

\end{document}